\documentclass[12pt]{article}
\usepackage{amsmath,amssymb}
\usepackage{graphicx}
\newcommand{\eqdef}{\stackrel{\text{def}}{=}}

\newcommand{\n}{\nonumber \\}
\newcommand{\bm}{\boldsymbol}
\newcommand{\ignore}[1]{}
\numberwithin{equation}{section}
\newcommand{\Romannumeral}[1]{\uppercase\expandafter{\romannumeral#1}}

\newtheorem{theo}{\bf Theorem}[section]

\newtheorem{rema}[theo]{\bf Remark}

\newtheorem{prop}[theo]{\bf Proposition}
\newtheorem{defi}[theo]{\bf Definition}

\allowdisplaybreaks[4]

\begin{document}

\baselineskip=20pt
\newcommand{\preprint}{
\vspace*{-20mm}\begin{flushleft}\end{flushleft}
}
\newcommand{\Title}[1]{{\baselineskip=26pt
  \begin{center} \Large \bf #1 \\ \ \\ \end{center}}}
\newcommand{\Author}{\begin{center}
  \large \bf 
  Ryu Sasaki${}$ \end{center}}
\newcommand{\Address}{\begin{center}
     Department of Physics and Astronomy, Tokyo University of Science,
     Noda 278-8510, Japan
        \end{center}}
\newcommand{\Accepted}[1]{\begin{center}
  {\large \sf #1}\\ \vspace{1mm}{\small \sf Accepted for Publication}
  \end{center}}

\preprint
\thispagestyle{empty}

\Title{Multivariate Meixner polynomials as  Birth and Death polynomials
}

\Author

\Address
\vspace{1cm}

\begin{abstract}
Based on the framework of Plamen Iliev, 
multivariate Meixner polynomials are constructed explicitly as Birth and Death polynomials.
They form the complete set of eigenpolynomials of a birth and death process 
with the birth and death rates at population 
$\bm{x}=(x_1,\ldots,x_n)\in\mathbb{N}_0^n$ are $B_j(\bm{x})=\bigl(\beta+\sum_{i=1}^nx_j\bigr)$ 
and $D_j(\bm{x})=c_j^{-1}x_j$,  $0<c_j$, $j=1,\ldots,n$, $\sum_{j=1}^nc_j<1$. 
The corresponding stationary distribution is  
$(\beta)_{\sum_{j=1}^nc_j}\prod_{j=1}^n(c_j^{x_j}/x_j!)(1-\sum_{j=1}^nc_j)^\beta$, 
the trivial $n$-variable generalisation of the orthogonality weight of the single variable Meixner polynomials.
The polynomials, depending on $n+1$ parameters ($\{c_i\}$ and $\beta$),
 satisfy the difference equation with the coefficients 
$B_j(\bm{x})$ and $D_j(\bm{x})$  $j=1,\ldots,n$,
which is the straightforward generalisation of the difference equation 
governing the single variable Meixner polynomials.
The polynomials  are  truncated $(n+1,2n+2)$ hypergeometric functions of Aomoto-Gelfand.
The polynomials and the derivation are very similar to those of the multivariate 
Krawtchouk polynomials reported recently.
\end{abstract}

%
%
\section{Introduction}
\label{sec:intro}
As the second member of the discrete multivariate hypergeometric orthogonal polynomials of Askey scheme
\cite{askey,ismail, koeswart,mKrawt}, {\em i.e.} those satisfying second order difference equations
with the nearest neighbour interactions,
the multivariate Meixner polynomials are constructed {\em explicitly} as
the  multivariate Birth and Death (BD) \cite{feller, KarMcG, ismail} polynomials.
Multivariate BD processes are the nearest neighbour interactions of multi-dimensional discrete systems,
the most basic type of interactions, like the well known Ising models.

The main strategy is essentially the same as that of the multivariate Krawtchouk 
orthogonal polynomials reported recently \cite{mKrawt}.
The basic framework, the orthogonality measure, the design of 
the truncated hypergeometric functions and the generating function connecting them,
is provided by a lucent paper of Iiiev \cite{I11} which follows the trends of many preceding 
works, {\em e.g.} \cite{Gri1,mizu}. The additional input is the reformulation of the birth and death
problems with the link to the associated self-adjoint matrices \cite{os12,bdsol,os34}.
It is markedly different from the traditional Karlin-McGregor approach \cite{KarMcG, ismail}.
Since the central motivation is the pursuit of multivariate orthogonal polynomials 
obeying second order difference equations, this paper is rather distant from other works 
on multivariate orthogonal polynomials
\cite{coo-hoa-rah77,tra,zheda,mizu1,mt,ilixu,IT,gr1,HR,khare,gr2,gr3,I,genest,xu,diaconis13,Gri2}.
 
 This paper is organised as follows.  Starting with a short r\'{e}sum\'{e} of the single variable
 Meixner polynomials in section \ref{sec:sMei}, the necessary ingredients for the construction of
 multivariate Meixner polynomials are listed in the logical order in section \ref{sec:nMeipath}.
In section \ref{sec:frame}, Iliev's framework for the multivariate Meixner polynomials 
is reproduced in my notation.
The main contents are in section \ref{sec:BDappr}.
An overview of multivariate birth and death problem is recapitulated in section \ref{sec:nBDview}.
The birth and death rates are introduced in section \ref{sec:nBDrates}.
The orthogonality of the Meixner polynomials is proved in section \ref{sec:nMort}.
The second order difference equation of the multivariate Meixner polynomials 
is demonstrated by using the generating function in section \ref{sec:Hteigen}.
The solution of the BD problem is presented in section \ref{sec:BDsol}.
The exceptional cases are remarked in section \ref{sec:excep}.
%
%
\section{Single variable  Meixner polynomials}
\label{sec:sMei}
Let us start with the summary of the single variable Meixner polynomials \cite{koeswart}
with real positive parameters $1>c>0$, $\beta>0$,
\begin{equation}
\qquad P_m(\beta,c;x)\eqdef{}_2F_1\Bigl(\genfrac{}{}{0pt}{}{-m,\,-x}{\beta}\Bigm|1-c^{-1}\Bigr)
=\sum_{k\in{\mathbb N}_0}\frac{(-m)_k(-x)_k}{(\beta)_k}\frac{(1-\frac1c)^{k}}{k!},
\quad m\in\mathbb{N}_0,
\label{sMeiform}
\end{equation}
which  are obtained by the generating function $G(\beta,c,x;t)$,
\begin{align} 
G(\beta,c,x;t)\eqdef\left(1-\frac{t}{c}\right)^x(1-t)^{-\beta-x}
=\sum_{m\in{\mathbb N}_0}\frac{(\beta)_m}{m!}P_m(\beta,c;x)t^n,
\label{sMeigen}
\end{align}
in which $(a)_n$ is the shifted factorial defined for $a\in\mathbb{C}$ and nonnegative integer $n$,
$(a)_0=1$, $(a)_n=\prod_{k=0}^{n-1}(a+k)$, $n\ge1$.
They satisfy the orthogonality relations with the normalised orthogonality weight $W(\beta,c;x)$,
\begin{align} 
 &\qquad W(\beta,c;x)\eqdef\frac{(\beta)_xc^x}{x!}(1-c)^{\beta},
 \label{sMeiW}\\
& \sum_{x\in{\mathbb N}_0}W(\beta,c;x)P_m(\beta,c;x)P_{m'}(\beta,c;x)
 =\frac{1}{(\beta)_n\frac{c^m}{m!}}\,\delta_{m\,m'},
\label{sMeiorth} 
\end{align}
and the second order difference equations are
\begin{align} 
 &\widetilde{\mathcal H}\eqdef(\beta+x)(1-e^\partial)+\frac{x}{c}(1-e^{-\partial}),
 \qquad \partial=\frac{d}{dx},\quad e^{\pm\partial}f(x)=f(x\pm1),
\label{sMeiHth}\\[2pt]
&\widetilde{\mathcal H}P_m(\beta,c;x)=\frac{1-c}{c}mP_m(\beta,c;x),\qquad m\in{\mathbb N}_0.
\label{sMeieq}
\end{align}
The eigenvalue is easily guessed as
\begin{equation*}
\widetilde{\mathcal H}x^m=m\left(-1+\frac{1}c\right)x^m+\text{lower degrees}, \qquad m\in{\mathbb N}_0.
\end{equation*}
It is straightforward to verify \eqref{sMeieq} in terms of the generating function,
\begin{equation}
\widetilde{\mathcal H}G(\beta,c,x;t)=\frac{1-c}{c}t\frac{\partial}{\partial t}G(\beta,c,x;t).
\label{htgensMei}
\end{equation}

%
%
\section{Path to multivariate Meixner polynomials}
\label{sec:nMeipath}
In this paper I present the explicit forms of multivariate Meixner polynomials satisfying second order difference equations
in a similar way to the single variable Meixner polynomials shown above.
The necessary ingredients are
\begin{enumerate}
\item (normalised) orthogonality weight $W(\bm x)$
\item generating function $G(\bm{x})$
\item general form of the polynomial $P_{\bm m}(\bm{x})$
\item proof of the orthogonality $\sum_{\bm{x}\in{\mathbb N}_0^n}W(\bm{x})P_{\bm m}(\bm{x})P_{{\bm m}'}(\bm{x})=0$, 
$\bm{m}\neq\bm{m}'$
\item second order difference operator $\widetilde{\mathcal H}$
\item eigenvalue spectrum $\mathcal{E}(\bm{m})$
\item proof of $\widetilde{\mathcal H}P_{\bm m}(\bm{x})=\mathcal{E}(\bm{m})P_{\bm m}(\bm{x})$
\end{enumerate}
Plamen Iliev provided the first four ingredients in \cite{I11}, which will be cited as I. 
The formulation of multivariate Birth and Death polynomials introduced by myself \cite{mKrawt}
provides the rest and the explicit forms of the multivariate Meixner polynomials are obtained.

%
%
\section{Framework for the multivariate Meixner polynomials due to Iliev \cite{I11}}
\label{sec:frame}
\begin{defi}
\label{IdefWG}
The normalised orthogonality weight of the multivariate Meixner polynomials is  a simple $n$-variable generalisation
of the single variable one \eqref{sMeiW}
\begin{equation}
W(\beta,{\bm c};{\bm x})\eqdef \frac{(\beta)_{|x|}\bm{c}^{\bm x}}{\bm{x}!}(1-|c|)^{\beta},
\quad \beta>0,
\label{nMeiWdef}
\end{equation}
in which the probability parameters $\{c_i>0\}$ are restricted by the summability of $W$,
\begin{align} 
{\bm x}&=(x_1,x_2,\ldots,x_n)\in{\mathbb N}_0^n,\quad |x|\eqdef\sum_{i=1}^nx_i,
\quad {\bm x}!\eqdef\prod_{i=1}^nx_i!,\n
{\bm c}&=(c_1,c_2,\ldots,c_n)\in{\mathbb R}_{>0}^n,\, \quad |c|\eqdef\sum_{i=1}^nc_i,
\quad {\bm c}^{\bm x}\eqdef \prod_{i=1}^nc_i^{x_i},\n
&\sum_{\bm{x}\in{\mathbb N}_0^n}W(\beta,{\bm c};{\bm x})=1
\quad \Rightarrow 0<|c|<1.
\end{align}
\end{defi}
First, prepare $n^2$ real parameters $\bm{u}=\{u_{i\,j}\}$, $i,j=1,\ldots,n$ constrained by the
conditions
\begin{align} 
&\sum_{i=1}^nc_iu_{i\,j}=|c|-1,\quad j=1,\ldots,n,
\label{ocond1}\\
&\sum_{i=1}^nc_iu_{i\,j}u_{i\,k}=|c|-1,\quad j\neq k,\quad j,k=1,\ldots,n.
\label{ocond2}
\end{align}
By introducing another closely related $n^2$ real parameters $\bm{b}=\{b_{i\,j}\}$, 
\begin{equation}
b_{i\,j}\eqdef 1-u_{i\,j},\qquad i,j=1,\ldots,n,
\label{budef}
\end{equation}
which are Iliev's original parameters, these conditions look simpler
\begin{align*} 
&\sum_{i=1}^nc_ib_{i\,j}=1,\quad j=1,\ldots,n,
\tag{I.2.3a}\\
&\sum_{i=1}^nc_ib_{i\,j}b_{i\,k}=1,\quad j\neq k,\quad j,k=1,\ldots,n.
\tag{I.2.3b}
\end{align*}
Similar conditions have been introduced by many authors for the construction of 
multivariate orthogonal polynomials \cite{Gri1,Gri11,mizu,Gri2,mKrawt}.
It should be stressed that these $n(n+1)/2$ conditions are not enough to determine
$n^2$ parameters ${\bm u}$ completely. 
\begin{defi}
\label{cbardef}
Iliev {\rm[I.(2.3b)]} also introduces another set of probability parameters
$\{\bar{c}_j\}$,  by
\begin{equation}
1-|c|+\sum_{i=1}^nc_iu_{i\,j}^2=\frac{1-|c|}{\bar{c}_j},\quad j=1,\ldots,n.
\label{barcdef}
\end{equation}
\end{defi}
\begin{defi}
\label{ngenfundef}
The generating function $G$ is defined by {\rm (I.(2.6))}
\begin{equation}
G(\beta,\bm{u},\bm{x};\bm{t})\eqdef (1-|t|)^{-\beta-|x|}\prod_{i=1}^n\left(1-\sum_{j=1}^nb_{i\,j}t_j\right)^{x_i},
\label{nGdef}
\end{equation}
in which
\begin{equation*}
{\bm t}=(t_1,t_2,\ldots,t_n)\in{\mathbb C}^n,\quad |t|\eqdef\sum_{i=1}^nt_i.
\end{equation*}
\end{defi}
\begin{defi}
\label{nmpoly}
The polynomials  $P_{\bm m}(\bm{x})$ are defined by the expansion of $G$ {\rm (I.(2.7))} around ${\bm t}={\bm 0}$,
\begin{equation}
G(\beta,\bm{u},\bm{x};\bm{t})
=\sum_{{\bm m}\in{\mathbb N}_0^n}\frac{(\beta)_{|m|}}{{\bm m}!}
P_{\bm m}(\beta,{\bm u};{\bm x}){\bm t}^{\bm m},\quad {\bm m}=(m_1,\ldots,m_n)\in{\mathbb N}_0^n,
\quad {\bm t}^{\bm m}=\prod_{i=1}^nt_i^{m_i}.
\label{GPolexp}
\end{equation}
\end{defi}
His two main Theorems are  on the orthogonality and the concrete form of $P_{\bm m}({\bm x})$.
\begin{theo}
\label{theoort}
{\bf Orthogonality} relation reads
\begin{align}
\sum_{{\bm x}\in{\mathbb N}_0^n}W(\beta,{\bm c};{\bm x})P_{\bm m}(\beta,{\bm u};{\bm x})
P_{{\bm m}'}(\beta,{\bm u};{\bm x})&=\frac{1}{\bar{W}(\beta,\bar{\bm c};{\bm m})}\,\delta_{{\bm m}\,{\bm m}'},
\quad {\bm m},{\bm m}'\in{\mathbb N}_0^n,
\label{ortrel}\\
\bar{W}(\beta,\bar{\bm c};{\bm m})&\eqdef \frac{(\beta)_{|m|}\bar{\bm c}^{\bm m}}{{\bm m}!}.
\label{barWdef}
\end{align}
\end{theo}
\begin{theo}
\label{nPform}
{\bf Concrete form of the polynomial in terms of ${\bf u}$} is 
\begin{gather}
\label{Pm}
P_{\bm{m}}(\beta,{\bm u};\bm{x})
\eqdef \sum_{\substack{\sum_{i,j}c_{ij}\\
(c_{ij})\in {\mathbb M}_{n}}}
\frac{\prod\limits_{i=1}^{n}(-x_{i})_{\sum\limits_{j=1}^{n}c_{ij}}
\prod\limits_{j=1}^{n}(-m_{j})_{\sum\limits_{i=1}^{n}c_{ij}}}
{(\beta)_{\sum_{i,j}c_{ij}}} \; \frac{\prod(u_{ij})^{c_{ij}}}{\prod c_{ij}!},
\end{gather}
in which ${\mathbb M}_{n}$ is the set of all $n\times n$ matrices with entries in ${\mathbb N}_0$.
This is the hypergeometric function of Aomoto-Gelfand {\rm \cite{AK,gelfand}} of type $(n+1,2n+2)$
and 
it  has almost  the same form as that of the multivariate Krawtchouk polynomials reported earlier
{\rm (\cite{mizu}.2)}, {\rm (\cite{diaconis13}.7)}, {\rm (\cite{mKrawt}.3.15)}, 
except that $-N$ is replaced by $\beta$.
\end{theo}

It should be stressed that these results constitute the framework only. In order to construct the 
multivariate Meixner polynomials explicitly, the $n\times n$ matrix of the 
 parameters ${\bm u}$ must be specified completely. 
 Two types of examples are presented in Iliev's article \cite{I11}\S5.

%
%
\section{Approach via Birth and Death problem}
\label{sec:BDappr}
 Now I present the procedure to derive the multivariate Meixner polynomials explicitly,
 {\em i.e.} to determine the parameters $\{u_{i\,j}\}$ explicitly,  based on the
multivariate Birth and Death problem setting.

%
%
\subsection{An overview of $n$-variate  Birth and Death problem}
\label{sec:nBDview}
This is a simple recapitulation of the problem setting of multivariate Birth and Death (BD) problems
reported in a previous paper \cite{mKrawt}.
For the single variable BD problem, the well-known approach by Karlin-McGregor \cite{KarMcG,ismail} is
based on the three term recursion relations of generic orthogonal polynomials.
My method \cite{bdsol} employs second order difference equations 
which are dual to the three term recursion relations.
All of the Askey scheme polynomials of a single discrete variable provide `exactly solvable BD problem' of one
population group. For the construction of multivariate orthogonal polynomials of Askey type,
I adopt the BD problem approach after the one successful example 
of the multivariate Krawtchouk polynomials \cite{mKrawt}.

The  differential equation for the birth and death (BD) process of  $n$ groups of unlimited population
 $\bm{x}=(x_1,x_2,\ldots,x_n)\in{\mathbb N}_0$ reads 
\begin{align}
\frac{\partial}{\partial t}\mathcal{P}(\bm{x};t)&=(L_{BD}\mathcal{P})(\bm{x};t)
=\sum_{\bm{y}\in{\mathbb N}_0^n}{L_{BD}}_{\bm{x}\,\bm{y}}\mathcal{P}(\bm{y};t),
\quad \mathcal{P}(\bm{x};t)\ge0,\quad \sum_{\bm{x}\in{\mathbb N}_0^n} \mathcal{P}(\bm{x};t)=1,
\label{bdeqformal2}\\
&=-\sum_{j=1}^n(B_j(\bm{x})+D_j(\bm{x}))\mathcal{P}(\bm{x};t)
+\sum_{j=1}^nB_j(\bm{x}-\bm{e}_j)\mathcal{P}(\bm{x}-\bm{e}_j;t)\n
&\hspace{5.8cm}+\sum_{j=1}^nD_j(\bm{x}+\bm{e}_j)\mathcal{P}(\bm{x}+\bm{e}_j;t),
\label{BDeq2}
\end{align}
in which $\bm{e}_j$ is the  $j$-th unit vector, $j=1,\ldots,n$. 
The birth and death rates for $n$ groups $B_j(\bm{x})$, $D_j(\bm{x})$  are all positive with the
boundary conditions
\begin{equation}
B_j(\bm{x})>0,\quad D_j(\bm{x})>0,\quad D_j(\bm{x})=0  \ \   \text{if}\ \ \bm{x}\in{\mathbb N}_0^n \ \ \text{and}\ \ \bm{x}-\bm{e}_j\notin{\mathbb N}_0^n    \quad j=1,\ldots,n.
\label{becond2}
\end{equation}
This is a typical example of the nearest neighbour interactions in $n$ dimensions.
The birth and death operator $L_{BD}$, an ${\mathbb N}_0^n\times{\mathbb N}_0^n$ matrix,  
can be expressed succinctly as
\begin{align}
L_{BD}&=-\sum_{j=1}^n\left[B_j(\bm{x})-B_j(\bm{x}-\bm{e}_j)e^{-\partial_j}+D_j(\bm{x})-D_j(\bm{x}+\bm{e}_j)e^{\partial_j}\right],
\label{LBDop0}\\
&=-\sum_{j=1}^n\left[(1-e^{-\partial_j})B_j(\bm{x})+(1-e^{\partial_j})D_j(\bm{x})\right],
\label{LBDop}\\
&=-\sum_{j=1}^n(1-e^{-\partial_j})\bigl(B_j(\bm{x})-D_j(\bm{x}+\bm{e}_j)\,e^{\partial_j}\bigr).
\label{LBDop1}
\end{align}
It is required that the system has a stationary distribution $W({\bm x})$
\begin{equation}
(L_{BD}W)(\bm{x})=0,\quad \sum_{\bm{x}\in{\mathbb N}_0^n}W(\bm{x})=1, \quad W(\bm{x})>0,
\quad \bm{x}\in{\mathbb N}_0^n,
\label{stationaryn}
\end{equation}
which constrains $\{B_j(\bm{x}),D_j(\bm{x})\}$ severely. 
The sufficient condition for the existence of the zero mode of $L_{BD}$ \eqref{stationaryn}, the stationary distribution $W(\bm{x})>0$, reads
\begin{align}
\bigl(B_j(\bm{x})-D_j(\bm{x}+\bm{e}_j)e^{\partial_j}\bigr)W(\bm{x})=0
\ \Rightarrow \frac{W(\bm{x}+\bm{e}_j)}{W(\bm{x})}=\frac{B_j(\bm{x})}{D_j(\bm{x}+\bm{e}_j)},
\quad j=1,\ldots,n,
\label{zerocond}
\end{align}
together with the compatibility conditions.
\begin{equation}
\frac{B_j(\bm{x})}{D_j(\bm{x}+\bm{e}_j)}\frac{B_k(\bm{x}+\bm{e}_j)}{D_k(\bm{x}+\bm{e}_j+\bm{e}_k)}
=\frac{B_k(\bm{x})}{D_k(\bm{x}+\bm{e}_k)}\frac{B_j(\bm{x}+\bm{e}_k)}{D_j(\bm{x}+\bm{e}_k+\bm{e}_j)},
\quad j,k=1,\ldots,n,
\label{compcond}
\end{equation}
When satisfied, these conditions determine entire $W(\bm{x})$ starting from the origin $W(\bm{0})$.
 
 By a similarity transformation of $L_{BD}$ in terms of $\sqrt{W({\bm x})}$, a new operator (matrix)
 $\mathcal H$ is introduced,
\begin{align}
\mathcal{H}&\eqdef -\bigl(\sqrt{W(\bm{x})}\bigr)^{-1}L_{BD}\sqrt{W(\bm{x})}.
\label{Hdef}\\
&=\sum_{j=1}^n\left[B_j(\bm{x})+D_j(\bm{x})-\sqrt{B_j(\bm{x})D_j(\bm{x}+\bm{e}_j)}\,e^{\partial_j}
-\sqrt{B_j(\bm{x}-\bm{e}_j)D_j(\bm{x})}\,e^{-\partial_j}\right],
\label{Hdef2}\\
\mathcal{H}_{\bm{x}\,\bm{y}}&=\sum_{j=1}^n\left[\bigl(B_j(\bm{x})+D_j(\bm{x})\bigr)\,\delta_{\bm{x}\,\bm{y}}-\sqrt{B_j(\bm{x})D_j(\bm{x}+\bm{e}_j)}\,\delta_{\bm{x}+\bm{e}_j\,\bm{y}}\right.\n
&\left.\hspace{5.2cm}
-\sqrt{B_j(\bm{x}-\bm{e}_j)D_j(\bm{x})}\,\delta_{\bm{x}-\bm{e}_j\,\bm{y}}\right].
\label{Hdef3}
\end{align} 
The operator $\mathcal{H}$ is a {\em positive semi-definite  real symmetric matrix} 
as is clear by the following factorisation,
\begin{align}
&\hspace{4cm}\mathcal{H}=\sum_{j=1}^n\mathcal{A}_j(\bm{x})^T\mathcal{A}_j(\bm{x}),
\qquad \mathcal{H}_{\bm{x}\,\bm{y}}=\mathcal{H}_{\bm{y}\,\bm{x}},
\label{Hfac}\\
&\mathcal{A}_j(\bm{x})\eqdef \sqrt{B_j(\bm{x})}-e^{\partial_j}\sqrt{D_j(\bm{x})},\ 
\mathcal{A}_j(\bm{x})^T= \sqrt{B_j(\bm{x})}-\sqrt{D_j(\bm{x})}\,e^{-\partial_j},\ 
j=1,\ldots,n.
\label{Ajdef}
\end{align} 
As $W({\bm x})$ is the zero mode of $L_{BD}$ \eqref{stationaryn}, 
$\sqrt{W(\bm{x})}$ is the zero mode of $\mathcal{A}_j(\bm{x})$ and $\mathcal{H}$
\begin{equation}
\mathcal{A}_j(\bm{x})\sqrt{W(\bm{x})}=0, \quad j=1,\ldots,n \quad \Longrightarrow \mathcal{H}\sqrt{W(\bm{x})}=0.
\label{zeromodes}
\end{equation}
Another operator $\widetilde{\mathcal H}$ is introduced by a similarity transformation of $\mathcal{H}$
in terms of the square root of the stationary distribution $W(\bm{x})$,
\begin{align}
\widetilde{\mathcal H}&\eqdef \bigl(\sqrt{W(\bm{x})}\bigr)^{-1}\mathcal{H}\sqrt{W(\bm{x})}.
\label{Hthdef}\\
&=\sum_{j=1}^n\left[B_j(\bm{x})\bigl(1-e^{\partial_j}\bigr)+D_j(\bm{x})\bigl(1-e^{-\partial_j}\bigr)\right],
\label{Hthdef2}
\end{align}
which provides the difference equations for the possible multivariate orthogonal polynomials of Askey type.
A trivial fact that a constant is the zero mode of $\widetilde{\mathcal H}$ is worth mentioning
\begin{equation}
\widetilde{\mathcal H}\,1=0,
\label{Hthzero}
\end{equation}
which corresponds to \eqref{zeromodes}.
%
%
\subsection{BD rates for Meixner}
\label{sec:nBDrates}
\begin{defi}{\bf BD rates for Meixner}
The following Birth and Death rates are adopted for the $n$-variate Meixner polynomials,
\begin{equation}
B_j({\bm x})\eqdef \beta+|x|,\qquad D_j({\bm x})\eqdef c_j^{-1}x_j,\quad c_j>0, 
\qquad j=1,\ldots,n,\quad \beta>0,
\label{BDMeidef}
\end{equation}
in which the parameters $\{c_j\}$, $j=1,\ldots,n$ are asuumed to be generic.
\end{defi}
They trivially satisfy the compatibility conditions \eqref{compcond} for $j,k=1,\ldots,n$,
\begin{align*}
\frac{B_j(\bm{x})}{D_j(\bm{x}+\bm{e}_j)}\frac{B_k(\bm{x}+\bm{e}_j)}{D_k(\bm{x}+\bm{e}_j+\bm{e}_k)}
=\frac{c_kc_k(\beta+|x|)(\beta+|x|+1)}{(x_j+1)(x_k+1)}
=\frac{B_k(\bm{x})}{D_k(\bm{x}+\bm{e}_k)}\frac{B_j(\bm{x}+\bm{e}_k)}{D_j(\bm{x}+\bm{e}_k+\bm{e}_j)}.
\end{align*}
They lead to the stationary distribution introduced in Definition \ref{IdefWG}
\begin{equation*}
W(\beta,{\bm c};{\bm x})\eqdef \frac{(\beta)_{|x|}\bm{c}^{\bm x}}{\bm{x}!}(1-|c|)^{\beta},
\quad \beta>0,
\tag{\ref{nMeiWdef}}
\end{equation*}
as the sufficient conditions \eqref{zerocond} are easily verified
\begin{align*}
\left((\beta+|x|)-c_j^{-1}(x_j+1)e^{\partial_j}\right)W(\beta,{\bm c};{\bm x})
=\frac{(\beta)_{|x|+1}\bm{c}^{\bm x}}{\bm{x}!}(1-|c|)^{\beta}-\frac{(\beta)_{|x|+1}\bm{c}^{\bm x}}{\bm{x}!}(1-|c|)^{\beta}=0.
\end{align*}
The summability of $W(\beta,{\bm c};{\bm x})$ due to the summation formula
\begin{equation}
\sum_{n=0}^\infty\frac{(\gamma)_n}{n!}z^n={}_1F_0\Bigl(\genfrac{}{}{0pt}{}{\gamma}{-}\Bigm|z\Bigr)
=(1-|z|)^{-\gamma}
\label{sumform}
\end{equation}
limits the parameter ranges
\begin{equation}
0<\sum_{j=1}^n c_j\equiv |c|<1.
\label{cjrange}
\end{equation}
\begin{defi}{\bf operator $\widetilde{\mathcal H}$ for $n$-variate Meixner} takes a very simple form
\begin{equation}
\widetilde{\mathcal H}=(\beta+|x|)\sum_{j=1}^n(1-e^{\partial_j})+\sum_{j=1}^nc_j^{-1}x_j(1-e^{-\partial_j}).
\label{HtMeix}
\end{equation}
\end{defi}
It is easy to see that the set of $n$-variate polynomials of maximal degree $M$ 
\begin{equation}
V_M({\bm x})\eqdef \text{Span}\{{\bm x}^{\bm m}|0\leq|m|\leq M\},\quad |m|\eqdef\sum_{j=1}^n m_j,
\end{equation}
is invariant under  $\widetilde{\mathcal H}$
\begin{equation}
\widetilde{\mathcal H}V_M({\bm x})\subseteq V_M({\bm x}),
\end{equation}
and  $\widetilde{\mathcal H}$ has eigenpolynomials in each $V_M({\bm x})$.
For generic values of the parameters $\{c_j\}$, these polynomials are orthogonal  with
each other due to the real symmetry of the operator $\mathcal{H}$ \eqref{Hfac}.

%
%
\subsection{Meixner polynomials are orthogonal with each other}
\label{sec:nMort}
It is easy to determine $n$ degree 1 eigenpolynomials of $\widetilde{\mathcal H}$ \eqref{HtMeix} with 
unknown coefficients $\{a_i\}$ and unit constant part,
\begin{align}
&P_{|m|=1}(\bm{x})=1+\sum_{i=1}^na_ix_i,
\quad
\widetilde{\mathcal H}P_{|m|=1}(\bm{x})=\lambda P_{|m|=1}(\bm{x}),
\label{deg1form}\\
&\Rightarrow -(\beta+|x|)\sum_{i=1}^na_i+\sum_{i=1}^nc_i^{-1}a_ix_i=\lambda\Bigl(1+\sum_{i=1}^na_ix_i\Bigr).
\nonumber
\end{align}
By equating the coefficients of $x_i$ and 1, an eigenvalue equations of $\{a_i\}$ are obtained,
\begin{align}
-\sum_{k=1}^na_k+c_i^{-1}a_i&=\lambda a_i,\qquad  -\beta\sum_{k=1}^na_k=\lambda,
\label{1eig}\\
&  \Longrightarrow \ a_i=\frac1\beta\frac{\lambda}{\lambda-c_i^{-1}}, \quad i=1,\ldots,n,.
\label{lambdaeq}
\end{align}
Here $\lambda$ is the  root of a degree $n$ characteristic polynomial $\mathcal{F}(\lambda)$
of an $n\times n$  matrix $F({\bm c})$  depending on $\{c_i\}$,
\begin{equation}
0=\mathcal{F}(\lambda)\eqdef \text{Det}\bigl(\lambda I_n-F({\bm c})\bigr),\quad  F({\bm c})_{i\,j}\eqdef -1+c_i^{-1}\delta_{i\,j}.
\label{chareq}
\end{equation}
For each eigenvalue $\lambda_j$, which is positive by construction,
the unknown coefficients $\{a_i\}$ are determined,
\begin{equation*}
a_{i,j}=\frac{\lambda_j}{\beta(\lambda_j-c_i^{-1})},\quad i,j=1,\ldots,n,
\end{equation*}
and it satisfies the relation
\begin{equation}
\sum_{i=1}^n\frac{1}{\lambda_j-c_i^{-1}}\equiv\sum_{i=1}^n\frac{c_i}{c_i\lambda_j-1}=-1,\quad j=1,\ldots,n.
\label{clambrel}
\end{equation}
Let us tentatively identify the above $j$-th solution as $\bm{m}=\bm{e}_j$ solution
\begin{equation}
P_{\bm{e}_j}(\bm{x})=1+\frac1\beta\sum_{i=1}^n\frac{\lambda_j}{\lambda_j-c_i^{-1}}x_i,\quad j=1,\ldots,n.
\label{e_jsol}
\end{equation}
By comparing these polynomials with the general hypergeometric functions \cite{AK, gelfand} in
{\bf Theorem \ref{nPform}}
\begin{gather*}
P_{\bm{m}}(\beta,{\bm u};\bm{x})
\eqdef \sum_{\substack{\sum_{i,j}c_{ij}\\
(c_{ij})\in {\mathbb M}_{n}}}
\frac{\prod\limits_{i=1}^{n}(-x_{i})_{\sum\limits_{j=1}^{n}c_{ij}}
\prod\limits_{j=1}^{n}(-m_{j})_{\sum\limits_{i=1}^{n}c_{ij}}}
{(\beta)_{\sum_{i,j}c_{ij}}} \; \frac{\prod(u_{ij})^{c_{ij}}}{\prod c_{ij}!},
\tag{\ref{Pm}}
\end{gather*}
the system parameters $\{u_{i\,j}\}$ are completely identified
\begin{equation}
P_{\bm{e}_j}(\beta,{\bm u};\bm{x})=1+\frac1\beta\sum_{i=1}^nu_{i\,j}x_i,\qquad 
u_{i\,j}=\frac{\lambda_j}{\lambda_j-c_i^{-1}},\quad i,j=1,\ldots,n.
\label{udef}
\end{equation}
The $n+1$ eigenvectors $\sqrt{W(\beta,{\bm c};{\bm x})}$,  
$\{\sqrt{W(\beta,{\bm c};{\bm x})}P_{\bm{e}_j}(\beta,{\bm u};\bm{x})\}$, $j=1,\ldots,n$ of the real 
symmetric matrix $\mathcal{H}$ \eqref{Hdef2} are orthogonal with each other for generic parameters $\{c_i\}$,
\begin{align}
\sum_{{\bm x}\in{\mathbb N}_0^n}W(\beta,{\bm c};{\bm x})P_{\bm{e}_j}(\beta,{\bm u};\bm{x})&=0,
\qquad \qquad \quad j=1,\ldots,n,
\label{0jort}\\
\sum_{{\bm x}\in{\mathbb N}_0^n}W(\beta,{\bm c};{\bm x})
P_{\bm{e}_j}(\beta,{\bm u};\bm{x})P_{\bm{e}_k}(\beta,{\bm u};\bm{x})&=0,\quad 
j\neq k, \quad j,k=1,\ldots,n.
\label{jkort}
\end{align}
By using the summation formulas
\begin{align*}
\sum_{{\bm x}\in{\mathbb N}_0^n}W(\beta,{\bm c};{\bm x})x_j&=\frac{c_j\beta}{1-|c|},
\qquad \qquad \qquad \qquad \quad \  j=1,\ldots,n,\\
\sum_{{\bm x}\in{\mathbb N}_0^n}W(\beta,{\bm c};{\bm x})x_jx_k&=\frac{\beta(\beta+1)c_jc_k}{(1-|c|)^2}+\frac{\beta c_j}{1-|c|}\,\delta_{j\,k},
\quad  j,k=1,\ldots,n,
\end{align*}
they are reduced to the same expressions as the constraining conditions \eqref{ocond1} and \eqref{ocond2}
of the $n^2$ real parameters $\bm{u}=\{u_{i\,j}\}$, $i,j=1,\ldots,n$  in {\bf Definition \ref{IdefWG}},
\begin{align*} 
&\sum_{i=1}^nc_iu_{i\,j}=|c|-1,\quad j=1,\ldots,n,
\tag{\ref{ocond1}}\\
&\sum_{i=1}^nc_iu_{i\,j}u_{i\,k}=|c|-1,\quad j\neq k,\quad j,k=1,\ldots,n.
\tag{\ref{ocond2}}
\end{align*}
The norm of $P_{\bm{e}_j}(\beta,{\bm u};\bm{x})$ can be calculated similarly
\begin{equation*}
\sum_{{\bm x}\in{\mathbb N}_0^n}W(\beta,{\bm c};{\bm x})P_{\bm{e}_j}(\beta,{\bm u};\bm{x})^2
=\frac1\beta\frac{1-|c|+\sum_{i=1}^nc_iu_{i\,j}^2}{1-|c|}=\frac1{\beta\bar{c}_j},
\qquad   j=1,\ldots,n.
\tag{\ref{barcdef}}
\end{equation*}

These lead to the following
\begin{theo}{\bf Orthogonality of $P_{\bm{m}}(\beta,{\bm u};\bm{x})$}\\
\label{Pmort}
The explicitly assembled parameters $\{u_{i\,j}=\lambda_j/(\lambda_j-c_i^{-1})\}$, $i,j=1,\ldots,n$ \eqref{udef}
satisfy all the constraints required for the definition of the generating function $G(\beta,\bm{u},\bm{x};\bm{t})$ 
 \eqref{nGdef} in {\bf Definition \ref{ngenfundef}}.
{\bf Definition \ref{nmpoly}} and {\bf Theorem \ref{theoort}, \ref{nPform}} assure that the polynomials 
$P_{\bm{m}}(\beta,{\bm u};\bm{x})$ \eqref{Pm} with the above  ${\bm u}=\{u_{i\,j}\}$
satisfy the orthogonality relation
\begin{align*}
\sum_{{\bm x}\in{\mathbb N}_0^n}W(\beta,{\bm c};{\bm x})P_{\bm m}(\beta,{\bm u};{\bm x})
P_{{\bm m}'}(\beta,{\bm u};{\bm x})&=\frac{1}{\bar{W}(\beta,\bar{\bm c};{\bm m})}\,\delta_{{\bm m}\,{\bm m}'},
\quad {\bm m},{\bm m}'\in{\mathbb N}_0^n,
\tag{\ref{ortrel}}\\
\bar{W}(\beta,\bar{\bm c};{\bm m})&= \frac{(\beta)_{|m|}\bar{\bm c}^{\bm m}}{{\bm m}!},
\tag{\ref{barWdef}}
\end{align*}
with the $\{\bar{c}_j>0\}$, $j=1,\ldots,n$ determined explicitly by \eqref{barcdef}.
\end{theo}

%
%
\subsection{All $\{P_{\bm m}(\beta,{\bm u};{\bm x})\}$ are  eigenvectors of $\widetilde{\mathcal H}$}
\label{sec:Hteigen}
Verifying  that all  higher degree ones $\{P_{\bm{m}}(\beta,{\bm u};\bm{x})\}$ \eqref{Pm} are also the eigenpolynomials of 
$\widetilde{\mathcal H}$ \eqref{HtMeix} is the next task.
The explicit forms of the eigenvalues $\mathcal{E}(\bm{m})$ 
are necessary. Since $P_{\bm m}(\beta,{\bm u};{\bm x})$ has the form
\begin{equation}
P_{\bm m}(\beta,{\bm u};\bm{x})=1+\frac1\beta\sum_{i,j=1}^nx_im_ju_{i\,j}+ \text{higher degrees},
\end{equation}
 $\widetilde{\mathcal H}$ acting on the higher degrees produces only the terms of linear and higher degrees.
The only constant part of $\widetilde{\mathcal H}P_{\bm{m}}(\beta,{\bm u};\bm{x})$ comes from 
$\beta\sum_{i=1}^n(1-e^{\partial_j})$ acting on the linear part,
\begin{equation*}
\beta\sum_{k=1}^n(1-e^{\partial_k})\left\{\frac1\beta\sum_{i,j=1}^nx_im_ju_{i\,j}\right\}
=-\sum_{i\,j}m_ju_{i\,j}=-\sum_{j=1}^nm_j\lambda_j\sum_{i=1}^n\frac{1}{\lambda_j-c_i^{-1}}
=\sum_{j=1}^nm_j\lambda_j,
\end{equation*}
in which \eqref{clambrel} is used.  
After applying $(1-e^{\partial_i})$, all higher degree terms vanish at the origin ${\bm x}={\bm 0}$ as they
consist of terms like $(x_i)_k(x_j)_l$, $k+l\ge2$, the typical structure of the hypergeometric functions.
This leads to the following
\begin{prop}
\label{lineareig}
If $P_{\bm{m}}(\beta,{\bm u};\bm{x})$ \eqref{Pm} is an eigenpolynomial of $\widetilde{\mathcal H}$ \eqref{HtMeix},
it has a linear spectrum
\begin{equation*}
\mathcal{E}(\bm{m})\eqdef\sum_{j=1}^nm_j\lambda_j.
\label{eigform}
\end{equation*}
\end{prop}

The remaining task is to demonstrate
\begin{equation}
\widetilde{\mathcal H}P_{\bm{m}}(\beta,{\bm u};\bm{x})
=\Bigl(\sum_{k=1}^nm_k\lambda_k\Bigr)P_{\bm{m}}(\beta,{\bm u};\bm{x}),\qquad 
\bm{x},\bm{m}\in{\mathbb N}_0^n.
\label{Pmeigpol}
\end{equation}
After the example of the single variable Meixner polynomials
in \S\ref{sec:sMei}, the following $n$-variables generalisation of the generating function formula \eqref{htgensMei} 
\begin{align}
\widetilde{\mathcal H}G(\beta,\bm{u},\bm{x};\bm{t})
&=\left(\sum_{k=1}^n\lambda_kt_k\frac{\partial}{\partial t_k}\right)G(\beta,\bm{u},\bm{x};\bm{t}),
\label{genform}\\
&G(\beta,\bm{u},\bm{x};\bm{t})\eqdef (1-|t|)^{-\beta-|x|}\prod_{i=1}^n\left(1-\sum_{j=1}^nb_{i\,j}t_j\right)^{x_i},
\tag{\ref{nGdef}}
\end{align}
will lead to \eqref{Pmeigpol} above.

The action of the two parts of $\widetilde{\mathcal H}$
\begin{equation*}
(\beta+|x|)\sum_{i=1}^n(1-e^{\partial_i}),\qquad
\sum_{i=1}^nc_i^{-1}x_i(1-e^{-\partial_i}),
\end{equation*}
on $G(u,\bm{x},t)$ is evaluated separately. The first part gives
\begin{align*}
&(\beta+|x|)\sum_{i=1}^n(1-e^{\partial_i})G(\beta,\bm{u},\bm{x};\bm{t})\\
 &\qquad =
(\beta+|x|)\sum_{i=1}^n\left(G(\beta,\bm{u},\bm{x};\bm{t})-\frac{1-\sum_{j=1}^nb_{i\,j}t_j}{1-|t|}G(\beta,\bm{u},\bm{x};\bm{t})\right)\\
&\qquad =(\beta+|x|)\frac1{1-|t|}\sum_{i,j=1}^n(b_{i\,j}-1)t_jG(\beta,\bm{u},\bm{x};\bm{t})\\
&\qquad = (\beta+|x|)\sum_{j=1}^n\frac{\lambda_jt_j}{1-|t|}G(\beta,\bm{u},\bm{x};\bm{t}), \hspace{45mm}   (*)
\end{align*}
in which \eqref{budef}, \eqref{udef} and \eqref{clambrel} are used to obtain
\begin{equation*}
\sum_{i=1}^n\sum_{j=1}^n(b_{i\,j}-1)t_j=-\sum_{i=1}^n\sum_{j=1}^nu_{i\,j}t_j
=-\sum_{i=1}^n\sum_{j=1}^n\frac{\lambda_j}{\lambda_j-c_i^{-1}}t_j=\sum_{j=1}^n\lambda_jt_j.
\end{equation*}
The second part gives
\begin{align*}
&\sum_{i=1}^nc_i^{-1}x_i(1-e^{\partial_i})G(\beta,\bm{u},\bm{x};\bm{t})\\
 &\qquad =\sum_{i=1}^nc_i^{-1}x_i\left(G(\beta,\bm{u},\bm{x};\bm{t})-\frac{1-|t|}{1-\sum_{j=1}^nb_{i\,j}t_j}G(\beta,\bm{u},\bm{x};\bm{t})\right)\\
 &\qquad =\sum_{i=1}^nc_i^{-1}x_i
 \frac{|t|-\sum_{j=1}^nb_{i\,j}t_j}{1-\sum_{j=1}^nb_{i\,j}t_j}G(\beta,\bm{u},\bm{x};\bm{t})\\
 &\qquad =\sum_{i,j=1}^n
 \frac{c_i^{-1}x_iu_{i\,j}t_j}{1-\sum_{j=1}^nb_{i\,j}t_j}G(\beta,\bm{u},\bm{x};\bm{t}).  \hspace{45mm} (**)
\end{align*}
Now r.h.s. of \eqref{genform} reads
\begin{align*}
(\beta+|x|)\sum_{k=1}^n\frac{\lambda_kt_k}{1- |t|}G(\beta,\bm{u},\bm{x};\bm{t})
 - \sum_{i,k=1}^n\frac{\lambda_kt_k x_ib_{i\,k}}{1-\sum_{j=1}^nb_{i\,j}t_j}G(\beta,\bm{u},\bm{x};\bm{t}).
\end{align*}
Here the following equality holds
\begin{equation*}
\lambda_kx_ib_{i\,k}=\lambda_kx_i\frac{-c_i^{-1}}{\lambda_k-c_i^{-1}}=-c_i^{-1}x_i\frac{\lambda_k}{\lambda_k-c_i^{-1}}=c_i^{-1}x_iu_{i\,k}.
\end{equation*}
The r.h.s. are qual to $(*)+(**)$ and this concludes the proof, which leads to the following
\begin{theo}
\label{eigproof}
The multivariate Meixner polynomials  $\{P_{\bm{m}}(\beta,{\bm u};\bm{x})\}$\! \eqref{Pm} constitute the complete 
set of eigenpolynomils of the difference operator $\widetilde{\mathcal H}$ \eqref{HtMeix} which is derived from
the $n$-variate Birth and Death process \eqref{BDMeidef}.
The orthogonality \eqref{ortrel} is the consequence of the self-adjointness of $\mathcal{H}$ \eqref{Hdef3} 
and the generic parameters.
\end{theo}
\begin{prop}{\bf $\mathfrak{S}_n$ Symmetry}
\label{Snsymm}
The Meixner polynomials  $\{P_{\bm{m}}(\beta,{\bm u};\bm{x})\}$\! \eqref{Pm} is 
invariant under the symmetric group $\mathfrak{S}_n$, 
due to the arbitrariness of the ordering of $n$ roots $\{\lambda_j\}$ 
of the characteristic equation \eqref{chareq} in the parameters $u_{i\,j}$ \eqref{udef}.
\end{prop}

%
%
\subsection{Solution of the BD problem}
\label{sec:BDsol}
In terms of the norm formula \eqref{ortrel}, the following  set of orthonormal vectors on the Hilbert space 
${\mathbb N}_0^n$ are
defined
\begin{align}
&\quad \sum_{\bm{x}\in{\mathbb N}_0^n}\hat{\phi}_{\bm{m}}(\beta,{\bm u};\bm{x})\hat{\phi}_{\bm{m}'}(\beta,{\bm u};\bm{x})
=\delta_{\bm{m}\,\bm{m}'}, \qquad \quad \bm{m},\bm{m}'\in{\mathbb N}_0^n,,
\label{ortrel2}\\
\hat{\phi}_{\bm{m}}(\beta,{\bm u};\bm{x})&\eqdef \sqrt{W(\beta,{\bm c};\bm{x})}P_{\bm{m}}(\beta,{\bm u};\bm{x})\sqrt{\bar{W}(\beta,\bar{\bm c};\bm{m})},
\qquad \bm{x},\bm{m}\in{\mathbb N}_0^n,
\label{ortphiderf}\\
&\bar{W}(\beta,\bar{\bm c};{\bm m})= \frac{(\beta)_{|m|}\bar{\bm c}^{\bm m}}{{\bm m}!},\quad
\sum_{\bm{m}\in{\mathbb N}_0^n}\bar{W}(\beta,\bar{\bm c};{\bm m})=(1-|\bar{c}|)^{-\beta}.
\label{baretaW}
\end{align}
They define an orthogonal matrix $\mathcal{S}$ on ${\mathbb N}_0^n$,
\begin{align}
&\hspace{4cm} \mathcal{S}_{\bm{x}\,\bm{m}}\eqdef \hat{\phi}_{\bm{m}}(\beta,{\bm u};\bm{x}),\n
&\bigl(\mathcal{S}^T\mathcal{S}\bigr)_{\bm{m}\,\bm{m}'}=
\sum_{\bm{x}\in{\mathbb N}_0^n}\bigl(\mathcal{S}^T\bigr)_{\bm{m}\,\bm{x}}\mathcal{S}_{\bm{x}\,\bm{m}'}
=\sum_{\bm{x}\in{\mathbb N}_0^n}\hat{\phi}_{\bm{m}}(\beta,{\bm u};\bm{x})\hat{\phi}_{\bm{m}'}(\beta,{\bm u};\bm{x})
=\delta_{\bm{m}\,\bm{m}'},\n
&\Rightarrow \delta_{\bm{x}\,\bm{y}}=\bigl(\mathcal{S}\mathcal{S}^T)_{\bm{x}\,\bm{y}}
=\sum_{\bm{m}\in{\mathbb N}_0^n}\mathcal{S}_{\bm{x}\,\bm{m}}\bigl(\mathcal{S}^T\bigr)_{\bm{m}\,\bm{y}}
=\sum_{\bm{m}\in{\mathbb N}_0^n}\hat{\phi}_{\bm{m}}(\beta,{\bm u};\bm{x})\hat{\phi}_{\bm{m}}(\beta,{\bm u};\bm{y}).
\label{dualort}
\end{align}
This means that $\hat{\phi}_{\bm{m}}(\beta,{\bm u};\bm{x})$ defines dual polynomials in $\bm{m}$ indexed by $\bm{x}$.
The discussion of dual polynomials is essentially the same as that for the multivariate
Krawtchouk polynomials in \cite{mKrawt} and it will not be repeated here.

In terms of the orthonormal eigenvectors $\{\hat{\phi}_{\bm{m}}(\beta,{\bm u};\bm{x})\}$ of $\mathcal{H}$ \eqref{Hdef2} and $L_{BD}$ \eqref{LBDop0} 
the transition probability  $\mathcal{T}({\bm x},{\bm y};t)$ of the Birth and Death equation \eqref{bdeqformal2} is expressed neatly.
The transition probability $\mathcal{T}({\bm x},{\bm y};t)$ is the solution of the BD equation \eqref{bdeqformal2} 
with the initial condition
\begin{equation*}
\mathcal{P}(\bm{x};0)=\delta_{\bm{x}\,\bm{y}},
\end{equation*}
\begin{align}
\mathcal{T}({\bm x},{\bm y};t)&=\hat{\phi}_{\bm{0}}(\beta,{\bm u};\bm{x})\hat{\phi}_{\bm{0}}(\beta,{\bm u};\bm{y})^{-1}
\sum_{{\bm m}\in{\mathbb N}_0^n}\,e^{-\mathcal{E}(\bm{m})t}\hat{\phi}_{\bm{m}}(\beta,{\bm u};\bm{x})\hat{\phi}_{\bm{m}}(\beta,{\bm u};\bm{y}),\quad t>0.
\label{Tran1}
\end{align}
It is trivial to verify the initial condition by \eqref{dualort}.  By using
\begin{equation*}
\hat{\phi}_{\bm{0}}(\beta,{\bm u};\bm{x})=\sqrt{W(\beta,{\bm c};\bm{x})}\sqrt{\bar{W}(\beta,\bar{\bm c};\bm{m})},
\end{equation*}
\eqref{Tran1} is reduced to
\begin{align*}
\mathcal{T}({\bm x},{\bm y};t)&=W(\beta,{\bm c};\bm{x})\bar{W}(\beta,\bar{\bm c};\bm{m})
\sum_{{\bm m}\in{\mathbb N}_0^n}\,e^{-\mathcal{E}(\bm{m})t}P_{\bm{m}}(\beta,{\bm u};\bm{x})P_{\bm{m}}(\beta,{\bm u};\bm{y}),\quad t>0,
\\
\frac{\partial}{\partial t}\mathcal{T}({\bm x},{\bm y};t)&=-W(\beta,{\bm c};\bm{x})\bar{W}(\beta,\bar{\bm c};\bm{m})
\sum_{{\bm m}\in{\mathbb N}_0^n}\,\mathcal{E}(\bm{m})\,e^{-\mathcal{E}(\bm{m})t}P_{\bm{m}}(\beta,{\bm u};\bm{x})P_{\bm{m}}(\beta,{\bm u};\bm{y}),\quad t>0.
\end{align*}
Since $P_{\bm{m}}(\beta,{\bm u};\bm{x})W(\beta,{\bm c};\bm{x})$ is the eigenvector of $L_{BD}$ with the
eigenvalue $-\mathcal{E}(\bm m)$,
\begin{equation*}
L_{BD}P_{\bm{m}}(\beta,{\bm u};\bm{x})W(\beta,{\bm c};\bm{x})=-\mathcal{E}(\bm m)
P_{\bm{m}}(\beta,{\bm u};\bm{x})W(\beta,{\bm c};\bm{x}),
\end{equation*}
the BD equation is satisfied
\begin{equation*}
\frac{\partial}{\partial t}\mathcal{T}({\bm x},{\bm y};t)=L_{BD}\mathcal{T}({\bm x},{\bm y};t),\quad t>0.
\end{equation*}
The transition probability \eqref{Tran1} ptrovides a simple example of Chapman-Kolmogorov equation in
${\mathbf N}_0^n$,
\begin{equation}
\mathcal{T}({\bm x},{\bm y};t+t')=\sum_{{\bm z}\in{\mathbb N}_0^n}\mathcal{T}({\bm x},{\bm z};t)
\mathcal{T}({\bm z},{\bm y};t').
\end{equation}
The r.h.s. is
\begin{align*}
&\sum_{{\bm z}\in{\mathbb N}_0^n}
\hat{\phi}_{\bm{0}}(\beta,{\bm u};\bm{x})\hat{\phi}_{\bm{0}}(\beta,{\bm u};\bm{z})^{-1}
\sum_{{\bm m}\in{\mathbb N}_0^n}\,e^{-\mathcal{E}(\bm{m})t}\hat{\phi}_{\bm{m}}(\beta,{\bm u};\bm{x})
\hat{\phi}_{\bm{m}}(\beta,{\bm u};\bm{z})\\
&\qquad \times 
\hat{\phi}_{\bm{0}}(\beta,{\bm u};\bm{z})\hat{\phi}_{\bm{0}}(\beta,{\bm u};\bm{y})^{-1}
\sum_{{\bm m}'\in{\mathbb N}_0^n}\,e^{-\mathcal{E}(\bm{m}')t'}\hat{\phi}_{\bm{m}'}(\beta,{\bm u};\bm{z})\hat{\phi}_{\bm{m}'}(\beta,{\bm u};\bm{y})\\
&=\sum_{{\bm z}\in{\mathbb N}_0^n}
\hat{\phi}_{\bm{0}}(\beta,{\bm u};\bm{x})\hat{\phi}_{\bm{0}}(\beta,{\bm u};\bm{y})^{-1}
\sum_{{\bm m}\in{\mathbb N}_0^n}\,e^{-\mathcal{E}(\bm{m})t}\hat{\phi}_{\bm{m}}(\beta,{\bm u};\bm{x})
\hat{\phi}_{\bm{m}}(\beta,{\bm u};\bm{z})\\
&\qquad \times 
\sum_{{\bm m}'\in{\mathbb N}_0^n}\,e^{-\mathcal{E}(\bm{m}')t'}\hat{\phi}_{\bm{m}'}(\beta,{\bm u};\bm{z})\hat{\phi}_{\bm{m}'}(\beta,{\bm u};\bm{y}).
\end{align*}
The   summation over ${\bm z}$, 
$\sum_{{\bm z}\in{\mathbb N}_0^n}\hat{\phi}_{\bm{m}}(\beta,{\bm u};\bm{z})
\hat{\phi}_{\bm{m}'}(\beta,{\bm u};\bm{z})=\delta_{{\bm m}\,{\bm m}'}$,
gives
\begin{align*}
r.h.s.&=
\hat{\phi}_{\bm{0}}(\beta,{\bm u};\bm{x})\hat{\phi}_{\bm{0}}(\beta,{\bm u};\bm{y})^{-1}
\sum_{{\bm m}\in{\mathbb N}_0^n}\,e^{-\mathcal{E}(\bm{m})(t+t')}\hat{\phi}_{\bm{m}}(\beta,{\bm u};\bm{x})
\hat{\phi}_{\bm{m}}(\beta,{\bm u};\bm{y})\\
&=l.h.s..
\end{align*}

%
%
\subsection{Exceptional Cases}
\label{sec:excep}
So far the parameter values are assumed to be generic. 
But obviously at certain parameter settings, 
the above hypergeometric formula \eqref{Pm} for the polynomials 
$\{P_{\bm{m}}(\beta,{\bm u};\bm{x})\}$ could go wrong. 
By construction, the eigenvalues $\{\lambda_i\}$ are positive 
and $\{c_i\}$ are also positive. 
Therefore, if the situation $\lambda_j=c_i^{-1}$  happens at some parameter setting, it leads to the
breakdown of the generic theory  as 
$u_{i\,j}=\frac{\lambda_j}{\lambda_j-c_i^{-1}}$  \eqref{udef} is ill-defined.

\subsubsection{$n=2$ Case}
\label{sec:n2exc}
The situation is most clearly seen when $n=2$. In this case the two eigenvalues are
the roots of the quadratic equation
\begin{align*}
&\lambda^2+(2-c_1^{-1}-c_2^{-1})\lambda+(1-c_1^{-1})(1-c_2^{-1})-1=0,\n
\lambda_1&=\frac12(-2+c_1^{-1}+c_2^{-1}-\Delta),\quad \lambda_2=\frac12(-2+c_1^{-1}+c_2^{-1}+\Delta),\\
& \Delta^2=4+(1/c_1-1/c_2)^2.
\end{align*}
When $c_1=c_2=c$, the eigenvalues are rational,
\begin{equation*}
c_1=c_2=c\ \Longrightarrow \lambda_1=-2+1/c>0,\quad \lambda_2=1/c,
\end{equation*}
and the singular situation occurs. That is, the general formula \eqref{Pm}  fails.
In this case, the degree 1 solution of $\widetilde{\mathcal H}$ \eqref{HtMeix}
corresponding to the eigenvalue $\lambda_2=1/c$ is
\begin{equation*}
P_{\bm{e}_1}(\bm{x})=const\times\bigl(x_1-x_2).
\end{equation*}
That is the constant part is vanishing and the assumption that degree one solutions have unit 
constant part \eqref{deg1form} simply fails. The situation is similar for general $n$ as 
stated by the following 
\begin{theo}
\label{theo:sing}
When some of the parameters $\{c_i\}$ coincide, the hypergeometric formula \eqref{Pm}
for the $n$-variate Meixner polynomials does not apply.
But the solutions of the difference equations $\widetilde{\mathcal H}$ \eqref{HtMeix}
still constitute the $n$-variate orthogonal polynomials.
\end{theo}
This is rather easy to see. If $c_j=c_k=c$, the matrix $c^{-1}I_n-F({\bm c})$ \eqref{chareq} 
has the $j$-th column $-(1,1,\ldots,1)^T$ and $k$-th column $-(1,1,\ldots,1)^T$,
thus the characteristic polynomial $\mathcal{F}(\lambda)$ vanishes at $\lambda=c^{-1}$.
In this case it is easy to show
\begin{equation*}
\widetilde{\mathcal H}(x_j-x_k)=\frac1c(x_j-x_k).
\end{equation*}
When  there exist $k$ identical $c_i$'s, $\mathcal{F}(\lambda)$ has a factor $(\lambda-c_i^{-1})^{k-1}$.
\begin{theo}
\label{theo:distinct}{\bf Distinct parameters $\{c_j\}$ are necessary}\\
All the parameters $\{c_j\}$ must be distinct for the hypergeometric formula \eqref{Pm}
for the $n$-variate Meixner polynomials to hold.
\end{theo}
\begin{rema}
\label{fullform}
It is a big challenge to derive a general formula of $n$-variate Meixner polynomials 
including all these exceptional cases.
\end{rema}


\end{document}